\documentclass[11pt]{article}
    \usepackage{amssymb, latexsym, amsmath}
    \newcommand{\setleftmargin}[1]{
      \addtolength{\textwidth}{\oddsidemargin}
      \addtolength{\textwidth}{1in}
      \addtolength{\textwidth}{-#1}
      \setlength{\oddsidemargin}{-1in}
      \addtolength{\oddsidemargin}{#1}
      \setlength{\evensidemargin}{\oddsidemargin}
    }
    \newcommand{\setrightmargin}[1]{
      \setlength{\textwidth}{8.5in}
      \addtolength{\textwidth}{-\oddsidemargin}
      \addtolength{\textwidth}{-1in}
      \addtolength{\textwidth}{-#1}
    }
    \newcommand{\settopmargin}[1]{
      \addtolength{\textheight}{\topmargin}
      \addtolength{\textheight}{1in}
      \addtolength{\textheight}{\headheight}
      \addtolength{\textheight}{\headsep}
      \addtolength{\textheight}{-#1}
      \setlength{\topmargin}{-1in}
      \addtolength{\topmargin}{-\headheight}
      \addtolength{\topmargin}{-\headsep}
      \addtolength{\topmargin}{#1}
    }
    \newcommand{\setbottommargin}[1]{
      \setlength{\textheight}{11in}
      \addtolength{\textheight}{-\topmargin}
      \addtolength{\textheight}{-1in}
      \addtolength{\textheight}{-\footskip}
      \addtolength{\textheight}{-#1}
    }
    \newcommand{\setallmargins}[1]{
      \settopmargin{#1}
      \setbottommargin{#1}
      \setleftmargin{#1}
      \setrightmargin{#1}
    } \setallmargins{1.2in}
\title{New Gelfond-Type Transcendental Numbers
}
\author{
Rachel M Chaphalkar, \thanks
{ Dept. of Math., Univ. of Wisconsin-Whitewater,
Whitewater, WI 53190, USA,  Email: chaphalr@uww.edu}
Suk-Geun Hwang, \thanks{
Department of Mathematics Education, Kyungpook University, Taegu 41566, Rep. of Korea and Department of Mathematics, Email:
sghwang@knu.ac.kr}
Choon Ho Lee, \thanks{
Department of Mathematics,
College of Natural Sciences, Hoseo University,
20 79th-Gil Hoseo-Ro Baebang-Eup,
Asan Choong-Nam,
336-795,
Korea, Email: chlee@hoseo.edu} and
Ki-Bong Nam \thanks{ Dept. of Math., Univ. of Wisconsin-Whitewater,
Whitewater, WI 53190, USA,  Email: namk@uww.edu: The fourth author was partially supported by Strategic Priorities Fund  of the Dept. of Math., UW-W. } }
\begin{document}
\maketitle
 \begin{abstract}
It is well known that value at a non-zero algebraic number of each of the functions $e^{x}, \ln x, \sin x, \cos x, \tan x, \csc x, \sec x, \cot x, \sinh x,$ $ \cosh x,$ $ \tanh x,$ and $\coth x$ 
is transcendental number (see Theorem 9.11 of \cite{N}). 
In this work, we show that let $f(x)$ be one of $e^{x}, \ln x, \sin x, \cos x, \tan x, \csc x, \sec x, $ $\cot x, $ $\sinh x,$ $ \cosh x,$ $ \tanh x,$ and $\coth x$ and
for a given equation $h_1(x)f(g(x))=h_2(x)$, if the equation has a non-zero solution $\alpha$ such that all of $f(g(\alpha))$, $h_1(\alpha),$ and $h_2(\alpha)$ are non-zeros, then the solution is a transcendental number where
$g(x),h_1(x), h_2(x)\in {\Bbb Q}[x]$, and $h_1(x)$ and $h_2(x)$ are relatively prime. We show that for given two transcendental numbers $\tau_1$ and $\tau_2,$ $\tau_1+\tau_2 i$ is a transcendental number. There are several methods to find transcendental numbers.
  \end{abstract}

      \newtheorem{lemma}{Lemma}
    \newtheorem{prop}{Proposition}
    \newtheorem{thm}{Theorem}
    \newtheorem{coro}{Corollary}
    \newtheorem{definition}{Definition}
\newtheorem{exa}{Example}[section]
\newtheorem{note}{Note}[section]


\section{Preliminaries}

Let ${\mathbb Z}$ be the set of all integers and ${\mathbb N}$ $(\hbox{resp}.$ ${\Bbb Z}^+)$ be the
set of all non-negative $(\hbox{resp. positive})$ integers. Let ${\mathbb Q}$ be the set of all rational numbers and ${\Bbb R}$ be the set of all real numbers.
%A function $f:{\Bbb R} \to {\Bbb R}$ is algebraic if for any algebraic number $\alpha$ of ${\Bbb R}$, $f(\alpha)$ is algebraic.
For a given transcendental numbers $\tau$ and given positive integers $n$ and $m$, we define the additive group
$$[\tau]_{n,m}=\{a_n^n\tau+\cdots +a_0+a_{-1}\tau^{-1}+\cdots +a_{-m}\tau^{-m}|a_n, \cdots , a_{-m}\in {\Bbb Q}\}.$$
Simillrily, we define
the subgroups $[\tau]_n=\{a_n^n\tau+\cdots +a_0|a_n, \cdots , a_{0}\in {\Bbb Q}\}$ and
$[\tau]^*_n=\{a_n^n\tau+\cdots a_1\tau|a_n, \cdots , a_{1}\in {\Bbb Q}\}$
of
$[\tau]_{n,m}$. 
Similarly we define a subgroup $[\tau]^*_{n,m}$ of $[ \tau]_{n,m}$.
As $n,m$ approach $\infty$, we define additive groups $[\tau]_{n\to \infty}=[\tau]_{\infty}$ and $[\tau]_{n,m\to \infty}=[\tau]_{\pm \infty}$.
The quadratic irrationality of the equation $ae^2+be+c=0$ in Theorem 1 is studied in the book \cite{Sh}.
Let $f(x)$ be one of $e^{x^j}, \sin x, \cos x, \tan x, \csc x, \sec x, \cot x, \sinh x, \cosh x, $ $\tanh x,$ and $\coth x$. If we define the commutative function ring ${\Bbb Q}[f(x)]=
\{a_nf(x)^n+\cdots +a_1f(x)+a_0|a_n, \cdots ,a_0\in {\Bbb Q}\}$, then 
${\Bbb Q}[f(x)]$ is isomorphic to the polynomial ring ${\Bbb Q}[x]$. 
We can define a degree on the mapping ring
${\Bbb Q}[f(x)]$ as the polynomial ring ${\Bbb Q}[x]$ (see \cite{L}).
The function ring ${\Bbb Q}[f(x)^{\pm 1}]$ is the Laurent extension of ${\Bbb Q}[f(x)]$. Similarly, we define a function ring with many variables.
For example, the zeros of the quadratic polynomials
$x^2+2x-3\in {\Bbb R}[x]$ and
$(\sin x)^2+2\sin x-3\in {\Bbb R}[\sin x]$ are totally different ones. 
The equations have different algebraic structures (see \cite{L}). 
But ${\Bbb R}[x]$ and ${\Bbb R}[\sin x]$ are formally isomorphic as ${\Bbb R}$-algebras.
Let ${\mathcal T}$ be the set of all the transcendental numbers and let us define $\sim$ on ${\mathcal T}.$
For given two transcendental numbers $\tau_1$ and $\tau_2$, we define $\sim$ as follows:\\
$\tau_1\sim \tau_2$ if there is an algebraic number $r$ such that $\tau_1=\tau_2+r$. Then $\sim$ is an equivalence relation (\cite{L1} and \cite{ST}). For a transcendental number
$\tau$, ${\overline \tau}$ denotes the equivalence class of $\tau$ with respect to $\sim$. Let $\overline {\mathcal T}$ be the set of all equivalence classes of ${\mathcal T}$ with respect to 
$\sim$.
A ring $R=\oplus _{a\in I} R_a$ is semi-graded if for any $a\in R_a$ and $b\in R_b$, $ab\in R_{a+b}\oplus R_{a-b}$ where $I$ is an additive index set, i.e., $I$ is a subset of ${\Bbb R}$.

%%%%%%%%%%%%

%%%%%%%%%

%\newpage

\section{New Transcendental Numbers}

%%%%%%%%%%%%%%%

%%%%%%%%%%

For any non-zero algebraic number $\alpha$, $e^{\alpha}, \ln \alpha, \sin \alpha, \cos \alpha, \tan \alpha, \csc \alpha, \sec \alpha, \cot \alpha, $ $\sinh \alpha, $ $\cosh \alpha, \tanh \alpha,$
and $\coth \alpha$ are transcendental numbers (see Theorem 9.11 of \cite{N}). In addition
$\ln \alpha$ is transcendental number where $\alpha$ is a non-zero algebraic number (see Theorem 9.11 of \cite{N}).
\begin{thm}
Let $f(x)$ be one of $e^{x}, \ln x, \sin x, \cos x, \tan x, \csc x, \sec x, \cot x, \sinh x,$ $ \cosh x,$ $ \tanh x,$ and $\coth x$. For $g\in {\Bbb Q}[f(x)]$ and $h(x)\in {\Bbb Q}[x]$,
if an equation $g=h(x)$ has a solution $\alpha$ such that $g(f(\alpha))\neq 0$ and $h(\alpha)\neq 0$, then $\alpha$ is a transcendental number.
\end{thm}
{\it Proof.}
Let $\alpha$ be a solution of the equation $g=h(x)$ in the theorem. If $\alpha$ is transcendental,
then there is nothing to prove. Assume that $\alpha$ is algebraic. Then $h(\alpha)$ is algebraic. However $g(f(\alpha))$ is transcendental. We have a contradiction.
Thus we have proven the theorem.
\quad $\Box$

\begin{thm}
Let $f(x)$ be one of $e^{x}, \ln x, \sin x, \cos x, \tan x, \csc x, \sec x, $ $\cot x, $ $\sinh x,$ $ \cosh x,$ $ \tanh x,$ and $\coth x$.
For a given equation $h_1(x)f(g(x))=h_2(x)$, if the equation has a non-zero solution $\alpha$ such that $f(g(\alpha))\neq 0$ and $h_1(\alpha)\neq 0\neq h_2(\alpha)$
with $h_1(x)$ and $h_2(x)$ are relatively prime, then the solution is a transcendental number where
$g(x),h_1(x), h_2(x)\in {\Bbb Q}[x]$. 
\end{thm}
{\it Proof.}
Since every complex number is either an algebraic number or a transcendental number, let us assume that a solution $\alpha$ of the equation of the theorem
is algebraic. Then 
$h_2(\alpha)$ is algebraic. Since $\alpha$ is algebraic, $h_1(\alpha)$ and $g(\alpha)$ are algebraic. However
$f(g(\alpha))$ is transcendental. This gives a contradiction. Thus $\alpha$ is transcendental.
\quad $\Box$\\

\begin{coro}
Let $\hbox {arc} \sin x=f(x),$ $arc \cos x=f(x),$ $arc \tan x=f(x),$ $arc \cot x=f(x),$ $arc \sec x=f(x),$ $arc \csc x=f(x)$ be given well-defined equations
where $f(x)\in {\Bbb Q}[x].$ If one of those equations has a non-zero solution $\alpha$ such that $f(\alpha)\neq 0$, then the solution is a transcendental number.
\end{coro}
{\it Proof.}
The proof of the corollary is straightforward by Theorem 1.
\quad $\Box$

\bigskip

\noindent {\bf Note 1.} \\
\noindent The real solution of $e=x-4$ is $x=e+4$ which is a transcendental number.
The real number solution of $e^x+x-12=0$ is a transcendental number, that is 
$2.27472787147\cdots.$ 
The equation $e^x+x-12=0$ has countably many complex number solutions. The equation $\pi^x+4x=49$ has a real solution
which is the transcendental number $x=3.14097\cdots $. The following is Thomas' example: the equation $2(\cos x)^2-\cos x-1=x^2-x$ has solutions $0$ and $0.4177913944\cdots$ which is a transcendental number
\cite{B}.
\quad $\Box$

\begin{prop}
Let $f(x)$ be one of $e^{x^j}, \sin x, \cos x, \tan x, \csc x, \sec x, \cot x, \sinh x, \cosh x,$ $ \tanh x,$ and $\coth x$ where $j$ is a fixed positive integer.
For a given equation $(f(x)+a_1)^k=g(x)$, if the equation has a real solution $\alpha$ such that $f(\alpha)\neq 0$, $f(\alpha)+a_1\neq 0,$ and $g(\alpha)\neq 0$, then the solution is a transcendental number where
$g(x)\in {\Bbb Q}[x]$, $a_1$ is a given algebraic number, and $k$ is a given positive integer. 
\end{prop}
{\it Proof.}
Let $\alpha$ be the solution of the equation of the proposition. Since $\sqrt [k] {\alpha}$ is algebraic, the proof of the proposition is straightforward by the proof of Theorem 1. 
\quad $\Box$

\begin{prop}
Let $f(x)$ be one of $e^{x^j}, \sin x, \cos x, \tan x, \csc x, \sec x, \cot x, \sinh x, \cosh x, \tanh x,$ and $\coth x$. For any $g\in {\Bbb Q}[f(x)]$, if $1\leq \deg (g)\leq 4$ and $h(x)\in {\Bbb Q}[x]$, if
an equation $g=h(x)$ has a non-zero solution $\alpha$ such that $g(f(\alpha))\neq 0$ and $h(\alpha)\neq 0$, then $\alpha$ is a transcendental number.
\end{prop}
{\it Proof.}
Without loss of generality, suppose $f(x)=\sin x$. Let $\alpha$ be a non-zero solution of the equation
$g(\sin x)=h(x)$ in the proposition. If $\alpha$ is a transcendental number, then there is nothing to prove.
If $f((\sin x))$ is linear, quadratic, cubic, or quartic, then
by the solution of a linear, the quadratic, Cardano's, or Ferarri's formulas. respectively, we have that the left side
$\sin \alpha$ is a transcendental number and the right side of the solution is algebraic. This contradiction shows that
$\alpha$ is a transcendental number.
\quad $\Box$

\bigskip

\noindent {\bf Note 2.} 
By Proposition 1, we know that for algebraic number $\sinh \alpha$ and $\cosh \alpha$ are transcendental numbers as well.
A real  number solution of $xe^x=-x+12$ is between $1.7$ and $1.8.$ The equation $e^x-x+7=0$ has no real solution, but it has countably many complex transcendental number solutions. One of the complex number solutions of
the equation $e^x-x+7=0$ is $x= 1.7701\cdots +  2.669613 \cdots \it i $ and its complex conjugate $x= 1.7701\cdots - 2.669613 \cdots \it i $ 
are transcendental solutions of the equation which has the minimal modulus of all the complex solutions of the equation. The complex solutions
of $e^x-x+7=0$ are discrete. The equation
$(e^{a_1x^j} +a_2)^k=f(x)$ has sometimes countably complex number solutions.
\quad $\Box$\\

\begin{prop}\label{tran}
Let $F$ be a field which is generated by all the algebraic numbers.
For any $f(x)\in F[x]$ and for any transcendental number $\tau$, $f(\tau)$ is a transcendental number where $f(x)\notin {\Bbb Q}[x]$.
Furthermore for any $f(x)\in F[x^{\pm 1}]$ and for any transcendental number $\tau$, $f(\tau)$ is a transcendental number.
\end{prop}
{\it Proof.}
For the proof of the proposition it is enough to prove for the Laurent polynomial $f(x)$ of 
$ F[x^{\pm 1}]$.
Let us put $f(x)=a_nx^n+\cdots +a_1x+a_0+a_{-1}x^{-1}+\cdots +a_{-m} x^{-m}$
where $a_n,\cdots ,a_{-m}\in F$, at least one of $a_n, \cdots ,a_1,a_{-1}, \cdots ,a_{-m}$ is not zero, and $n,\cdots, m$ are positive integers or zeros. We have that
$f(\tau)=a_n\tau^n+\cdots +a_1\tau+a_0+a_{-1}\tau^{-1}+\cdots +a_{-m} \tau^{-m}$
and set $f(\tau)$ equals to $\alpha$ that is $\alpha=a_n\tau^n+\cdots +a_1\tau+a_0+a_{-1}\tau^{-1}+\cdots +a_{-m} \tau^{-m}.$ If $\alpha$ is a transcendental number,
then there is nothing to prove. Let us assume that $\alpha$ is algebraic. We have that
$a_n\tau^{n+m}+\cdots +a_1\tau^{m+1}+a_0\tau^m+a_{-1}\tau^{m-1}+\cdots +a_{-m} =\alpha \tau^m.$ Since $a_n,\cdots ,a_{-m}, \alpha$ are algebraic.
This implies that $\tau$ is algebraic. We have a contradiction. Thus $f(\tau)$ is a transcendental number.
\quad $\Box$

\begin{coro}\label{Fer}
If $\tau$ is a transcendental number and $\alpha_0,\cdots, \alpha_4$ are algebraic numbers, then
$f_1(\tau)=\alpha_4 \tau^{4}+\alpha_3\tau^{3}+\alpha_2\tau^{2}+\alpha_1 \tau+\alpha_0$,
$f_2(\tau)=\alpha_4 \tau^{3}+\alpha_3\tau^{2}+\alpha_2\tau+\alpha_1 \tau^{-1}+\alpha_0$,
$f_3(\tau)=\alpha_4 \tau^{2}+\alpha_3\tau^{1}+\alpha_2\tau^{-1}+\alpha_1 \tau^{-2}+\alpha_0$,
$f_4(\tau)=\alpha_4 \tau+\alpha_3\tau^{-1}+\alpha_2\tau^{-2}+\alpha_1 \tau^{-3}+\alpha_0$,
and
$\alpha_4 \tau^{-1}+\alpha_3\tau^{-2}+\alpha_2\tau_3^{-3}+\alpha_1 \tau^{-4}+\alpha_0$
are transcendental numbers where at least one of $\alpha_1,\cdots, \alpha_4$ is not zero.
\end{coro}
{\it Proof.}
The proof of the corollary is straightforward by Proposition 3 and the another proof of the corollary is the following.
If $\beta$ is the number of one of the equations of the proposition, then
since $\alpha_0, \alpha_1,\alpha_2,\alpha_3, \alpha_4$ are algebraic, by Ferrari's quartic formula, a solution of the equation $f_i(\tau),$ $1\leq i\leq 4,$ is algebraic.
But $\tau$ is transcendental, thus we have a contradiction. So $\beta$ is a transcendental number.
\quad $\Box$\\

\begin{coro}
If  $a,b,c,$ and $\alpha$ are given algebraic numbers, then $a(\cos \alpha)^2+b\sin \alpha+c$, $a\sin \alpha ^2+b\cos \alpha+c$,
$a(\tan \alpha)^2+b\sec \alpha+c$, and $a(\cot \alpha)^2+b\csc \alpha+c$ are transcendental numbers where at least one of $a$ and $b$ is not zero, and
$a(\cos \alpha)^2+b\sin \alpha+c$, $a\sin \alpha ^2+b\cos \alpha+c$,
$a(\tan \alpha)^2+b\sec \alpha +c$, and $a(\cot \alpha)^2+b\csc \alpha +c$ are not zeros.
\end{coro}
{\it Proof.}
The proof of the corollary is straightforward by Proposition 3 and the elementary trigonometric identities.
\quad $\Box$\\

%\begin{coro}
%Every element $l$ of the rings ${\bf <}\tau {\bf >}_{ n}$ and ${\bf <}\tau{\bf >}_{n,m}$
%is a transcendental number where $\tau$ is a given transcendental number and $l\notin {\Bbb Q}$.
%\end{coro}
%{\it Proof.}
%The proof of the corollary is straightforward by the proofs of Proposition 3.
%\quad $\Box$

\begin{coro}
For any element $l$ of the group $[\tau]_{\pm \infty},$ $l$ is a transcendental number where $l\notin {\Bbb Q}$ and $\tau$ is a given transcendental number.
For any non-zero element $l$ of the group $[\tau]_{\pm \infty}^*,$ $l$ is a transcendental number.
\end{coro}
{\it Proof.}
The proof of the corollary is straightforward by Proposition 3.
\quad $\Box$

\begin{prop}
For any two different integers $i$ and $j,$ and a given transcendental number $\tau$, $\overline {\tau^i} \cap \overline {\tau^j}=\phi$ and $\cup_{k\in {\Bbb N}} \overline {\tau^k} \subset {\bf <}\tau {\bf >}_{n,m}$ hold.
The cardinality of $\overline {\mathcal T}$ is uncountable. Furthermore, the set $\overline {\mathcal \tau}=\{\tau+\alpha | \alpha \hbox { is an algebraic number}\}$ is dense in ${\Bbb R}.$
\end{prop}
{\it Proof.} 
Let $\tau$ be the transcendental number of the proposition.
Since $\tau^i-\tau^j$ is a transcendental number, $\overline {\tau^i}$ and $ \overline {\tau^j}$ are disjoint. For any element $l$ of ${\overline \tau^k}$, $l=\tau^k+r$ where $r$ is an algebraic number.
This implies that $l$ is an element of ${\bf <}\tau {\bf >}_{n,m}.$
For the second statement of the proposition, for any $\tau_1\in \overline {\tau}\in \overline {\mathcal T}$, we have that $\tau_1=\tau +\alpha$ for an algebraic number $\alpha$. Since $ {\mathcal T}$ 
is uncountable and $ \overline {\tau}$ is countable, $\overline {\mathcal T}$ is uncountable. For the density of $\overline {\tau}$, for any given $\tau_1$ of 
$\overline {\tau}$ and given $\epsilon>0,$ there is
a positive integer $n$ such that $\frac{1}{n}<\epsilon$. This implies that the transcendental number $\tau_1+\frac{1}{n}$ of $\overline {\tau}$ is an element $(\tau -\epsilon, \tau+\epsilon).$ This implies that
$\overline {\tau}$ is dense in ${\Bbb R}$.
This completes the proof of the proposition.
\quad $\Box$

\begin{prop}
For any given transcendental numbers $\tau_1$ and $\tau_2$, $\tau_1+\tau_2i$ and $\tau_1-\tau_2i$ are transcendental numbers.
Specifically, $e+\pi i$ and $e-\pi i$ are transcendental numbers.
\end{prop}
{\it Proof.}
Symmetrically let us assume that $\tau_1+\tau_2i$ is an algebraic number. This implies that $\tau_1-\tau_2i$ is algebraic.
Assume these algebraic numbers to be $\alpha_1$ and $\alpha_2$ respectively.
By adding them, we have that $2\tau_1=\alpha_1+\alpha_2.$ This implies that $2\tau_1$ is a transcendental number but
$\alpha_1+\alpha_2$ is an algebraic number. We have a contradiction. This implies that
$\tau_1+\tau_2i$ is a transcendental number. 
Symmetrically $\tau_1-\tau_2i$ is a transcendental number.
This completes the proofs of the proposition.
\quad $\Box$\\

\begin{coro}
All most complex numbers of the complex number plane ${\Bbb C}$ are transcendental numbers as ${\Bbb R}$.
\end{coro}
{\it Proof.}
The proof of the corollary is straightforward by Proposition 5.
\quad $\Box$

\bigskip

\noindent As a result of Proposition 5, for any given real transcendental numbers $\tau_1$ and $\tau_2$, at least one of $\tau_1+\tau_2$ or $\tau_1-\tau_2$ is a transcendental number.
It is an interesting problem for given non-zero integers $i$ and $j$, whether both of $\pi^i+e^j$ and $\pi^i-e^j$ are transcendental numbers or not.
%For a given transcendental number $\tau$, $[\tau]_{\pm \infty}$ is ${\Bbb N}$-graded as follows:
%$${\bf <} \tau{\bf >}_{\pm \infty}=\bigoplus _{n\in {\Bbb N}} {\bf <} \tau {\bf >}_{n,n}$$

\begin{prop}
For a given transcendental number, the additive groups $[\tau]_n$, $[\tau]_n^*$, $[\tau]_{n,m}$, $[\tau]_{n,m}^*$, 
$[\tau]_{\infty}$, and  $[\tau]_{\pm \infty}$ are
${\Bbb Z}$-modules and they are ${\Bbb Q}$-modules. 
%The ring $<\tau>_{\pm \infty}$ is semi-graded.
\end{prop}

\noindent We have the following obvious results.

\begin{coro}
For a given transcendental number $\tau$, if $n_1\neq n_2$ or $m_1\neq m_2,$ then the vector spaces $[ \tau]_{n_1,m_1}$ and
$[ \tau]_{n_1,m_1}$ over ${\Bbb Q}$
are not isomorphic. For every element $l$ of $[ \tau]_{n_1,m_1}$, $l$ is a transcendental number where $l\notin {\Bbb Q}.$
\end{coro}
{\it Proof.}
The proofs of the corollary is straightforward by Proposition 3.
\quad $\Box$

\begin{prop}
For a given transcendental number $\tau$, the rings ${\bf <}\tau {\bf >}_{\infty}$
$($resp. ${\bf <}\tau {\bf >}_{\pm \infty})$ and the polynomial ring ${\Bbb Q}[x]$
$($resp. the Laurent extension ${\Bbb Q}[x^{\pm 1}]$ of the polynomial ring ${\Bbb Q}[x])$ are isomorphic.
Consequently the polynomial ring ${\Bbb Q}[x]$ and the Laurent extension ${\Bbb Q}[x^{\pm 1}]$ are integral domains.
\end{prop}
{\it Proof.}
The proof of the proposition is obvious, so let us omit it.
\quad $\Box$\\

\bigskip

\noindent {\bf Open question:}\\
For given transcendental numbers $\tau_i$, $1\leq i\leq 4,$ is the number 
$\tau_1+\tau_1 {\mathbf i}+ \tau_2 {\mathbf j}+\tau_3 {\mathbf k}$ a transcendental number in the 
quaternion ${\mathbf Q}$ (see \cite{E})?\\

\bigskip

\noindent {\bf Acknowledgement:}
The Authors thank you Thomas Preu for his reviewing in the theorems and corollaries of the manuscript
and his interesting examples for improving the manuscript.


\begin{thebibliography}
\frenchspacing
\parskip 11pt


\bibitem{B} A. Baker, {\it Transcendental Number
Theory}, Springer-Verlag, 1982.

\bibitem{E} Samuel Eilenberg and Ivan Niven, {\it The fundamental theorem of algebra for quaternions}, 
Bull. Amer. Math. Soc. 50(4): 246-248 (April 1944).\\

\bibitem{L} Jongwoo Lee, Xueqing Chen, Seul Hee Choi, and
Ki-Bong Nam, {\it Automorphism groups of some stable Lie algebras,}
Journal of Lie Theory, Heldermann Verlag, Vol. 21, 2011, 457-468.\\

\bibitem{L1} William J. LeVeque, {\it Fundamentals of Number Theory},
Cambridge University Press, London, 1975.

\bibitem{N} Ivan Niven, {\it Irrational Numbers}, The Carus Mathematical Monographs,
Number 12, 1967.\\

\bibitem{P} Thomas Preu, {\it Personal Communications}, 2021.

\bibitem{Sh} A. B. Shidlovskii, {\it Transcendental Numbers}, Translated from Russian by N. Koblitz, Walter de Gruyter, Berlin New York, 1989.


\bibitem{ST} Ian Stewart, {\it Algebraic Number Theory},
Chapman and Hall, 1979.

\end{thebibliography}
\end{document}